\input amstex
\documentstyle{amsppt}
\magnification=\magstep1 \NoRunningHeads
\topmatter

\title
 Furstenberg entropy~values~for nonsingular actions of~groups without property (T)
 \endtitle

\author
Alexandre I. Danilenko
\endauthor

\email
alexandre.danilenko@gmail.com
\endemail

\address
 Institute for Low Temperature Physics
\& Engineering of National Academy of Sciences of Ukraine, 47 Lenin Ave.,
 Kharkov, 61164, UKRAINE
\endaddress
\email alexandre.danilenko\@gmail.com
\endemail

\abstract
Let $G$ be a discrete countable infinite group that does not have Kazhdan's property ~(T) and let $\kappa$ be a  generating probability measure on $G$.
Then  for each $t>0$, there is a  type $III_1$ ergodic free nonsingular $G$-action  whose   $\kappa$-entropy (or the Furstenberg entropy) is $t$.
%If $W$ is trivial, the same is true but the type of  the action is $II_\infty$.
\endabstract

 \loadbold

\NoBlackBoxes

\endtopmatter
\document

\head 0. Introduction
\endhead
Let $G$ be a discrete countable infinite group.
 A probability measure $\kappa$ on $G$ is called {\it generating} if the support  of $\kappa$ generates $G$ as a semigroup.
 Let $T=(T_g)_{g\in G}$ be a nonsingular action of $G$ on a standard probability space $(X,\goth B,\mu)$.
 The {\it Furstenberg entropy}    (or $\kappa$-{\it entropy}) of $T$ is defined by
$$
h_\kappa(T,\mu): = -\sum_{g\in G}\kappa(g)\int_X\log\frac{d\mu\circ T_g}{d\mu}(x)d\mu(x)
 $$
  (see  \cite{Fu}).
%We will write $h_\kappa(T,\mu)$ in place of $h_\kappa(T)$ if we need to show which quasiinvariant measure is under question.
 Jensen's inequality implies that $h_\kappa(T,\mu)\ge 0$ and that (for generating measures) equality holds if and only if $\mu$ is  invariant under $T$.
 Of course, the $\kappa$-entropy is invariant under conjugacy.
 If $\sum_{g\in G}\kappa(g)\frac{d\mu\circ T_g}{d\mu}(x)=1$ for a.e. $x\in X$ then $T$ is called {\it $\kappa$-stationary}.
 {\it Furstenberg entropy realization problem} is to describe all values that $\kappa$-entropy takes on the set of $\kappa$-stationary actions.
 The problem appears  quite difficult.
 Some progress was achieved in recent papers \cite{NeZi}, \cite{Ne}, \cite{Bo}, \cite{HaTa}.
 To state one of the results  on the entropy realization problem we first recall that 
 $G$ has Kazhdan's {\it property (T)} if every unitary representation of $G$ which has almost invariant vectors admits a nonzero invariant vector (see \cite{Be--Va}).
 It was shown in \cite{Ne} that if $G$ has property (T) then for every  generating measure $\kappa$, the pair $(G,\kappa)$
 has an {\it entropy gap}, i.e. there 
 exists some constant $\epsilon=\epsilon(G,\kappa)>0$ such that the $h_\kappa(T)>\epsilon$ for each purely infinite ergodic   stationary $G$-action $T$.
 We recall that an ergodic action is called {\it purely infinite} if
  it does not admit an equivalent invariant probability measure.
In \cite{Bo--Ta} the converse statement was proved: if $T$ does not have property (T) then
for each generating measure $\kappa$,
$$
\inf\{h_\kappa(T,\mu)\mid T\text{ is purely infinite, ergodic, $\mu$-nonsingular action of $G$}\}=0.\tag0-1
$$
We note that the authors of \cite{Bo--Ta} consider $\kappa$-entropy values on arbitrary (not only stationary as in the other aforementioned papers) purely infinite nonsingular actions.
They also show that the entropy gap for $(G,\kappa)$ established in \cite{Ne} for the stationary actions holds also for all (purely infinite) ergodic nonsingular actions.
In this connection we note that if a purely infinite  ergodic  $G$-action is stationary  then the space of the action is non-atomic.
However, in the general (non-stationary) case considered in \cite{Bo--Ta}, there exist purely infinite  transitive $G$-actions on purely atomic measure spaces.
 In particular, the action of $G$ on itself via rotations is free, nonconservative, ergodic and purely infinite.
 We consider such actions as {\it pathological}.
 Unfortunately, the proof the main result from \cite{Bo--Ta} does not exclude appearance  of pathological actions in \thetag{0-1}.

Our purpose  in the present paper is to refine the main result from \cite{Bo--Ta} in two aspects: to examine all possible values for the $\kappa$-entropy and  ``get rid'' of possible pathological actions on which such values are attained.
In fact, we show more.
\comment
the  $\kappa$entropy takes {\it all possible} positive values.
 Secondly, we prove that this happens on the subclass of  ergodic $G$-actions  {\it of type $III_1$ }as well as on the subclass of ergodic $G$-actions {\it  of type $III_\lambda$, $0<\lambda<1$}.
 \endcomment

\proclaim{Main Theorem}
  Let $G$ do not have property (T).
  Let $\kappa$ be a generating measure on $G$.
  Then the following are satisfied.
  \roster
  \item
For  each real $t\in(0,+\infty)$, there is a  type $III_1$ ergodic free nonsingular action $T=(T_g)_{g\in G}$ on a standard probability space $(X,\mu)$  such that   $h_\kappa(T,\mu)=t$.
  \item
  For  each real $t\in(0,+\infty)$, there is  $\lambda\in(0,1)$ and a  type $III_\lambda$ ergodic free nonsingular action $T=(T_g)_{g\in G}$ on a standard probability space $(X,\mu)$  such that
   $h_\kappa(T,\mu)=t$.
  \endroster

\comment
  \roster
  \item
  For each $\epsilon>0$ and $\lambda\in [0,1]$, there is an ergodic free nonsingular  action $T$ of  $G$ of  Krieger's type $III_\lambda$ such that $0<h_\kappa(T)<\epsilon$.
  \item
  More generally,
  let $T$ be an ergodic but not strictly ergodic free nonsingular  action of $G$ on a non-atomic probability space.
  Then for each $\epsilon>0$ and, there is a type $III$ ergodic nonsingular extension $\widehat T$ of  $T$ such that  $h_\kappa(T)<h_\kappa(\widehat T)<h_\kappa(T)+\epsilon$.
  \endroster
  \endcomment
 \endproclaim

The proof is based on the measurable orbit theory (see \cite{FeMo}, \cite{Sc1} and a survey \cite{DaSi}) and cohomology properties of non-strongly ergodic actions \cite{Sc2},  \cite{JoSc}.

\subsubhead Acknowledgements
\endsubsubhead
I thanks K. Schmidt and S. Sinelshchikov for their useful remarks that I used in writing Appendix~A.

\head 1. Some background on orbit theory
\endhead
Let $T$ be an ergodic free nonsingular action of $G$ on a standard nonatomic probability space $(X,\goth B,\mu)$.
 Denote by $\Cal R$ the $T$-orbit equivalence relation on $X$.
 We recall that the {\it full group} $[\Cal R]$ of $\Cal R$ consists of all one-to-one nonsingular transformations $r$ of $(X,\mu)$ such that the graph of $r$ is  a subset of  $\Cal R$.
 Given a locally compact second countable  group  $H$, denote by $\lambda_H$ a left Haar measure on $H$.
 A Borel map $\alpha:\Cal R\to H$ is called a {\it cocycle} of $\Cal R$ if $\alpha(x,y)=\alpha(x,z)\alpha(z,y)$ for all points  $x,y,z$ from a $\mu$-conull subset of $X$ such that $x\sim_\Cal R y\sim_\Cal Rz$.
 By $T(\alpha)=(T(\alpha)_g)_{g\in G}$ we denote the $\alpha$-{\it skew product extension of $T$}, i.e.
 a $G$-action on the product space $(X\times H,\mu\times\lambda_H)$:
 $$
 T(\alpha)_g(x,h)=(T_gx,\alpha(T_gx,x)h).
 $$
 It is obvious that $T(\alpha)$ is $(\mu\times\lambda_H)$-nonsingular.
We say that $\alpha$ is {\it ergodic} if $T(\alpha)$ is ergodic.
Consider the $H$-action on $(X\times H,\mu\times\lambda_H)$ by rotations (from the right) along the second coordinate.
It commutes with $T(\alpha)$.
The restriction of this action to the sub-$\sigma$-algebra
 of $T(\alpha)$-invariant Borel subsets is called {\it the action of $H$ associated  with $\alpha$}. 
 It is ergodic.
 It is trivial if and only if $\alpha$ is ergodic.
 If $H=\Bbb R^*_+$ and $\alpha(T_gx,x):=\log\frac{d\mu\circ T_g}{d\mu}(x)$ at a.e. $x$ for each $g\in G$ then $\alpha$ is called {\it the Radon-Nikodym cocycle} of $\Cal R$.
 It does not depend on the choice of nonsingular group action generating $\Cal R$.
 The corresponding associated action of $\Bbb R_+^*$ is called {\it the associated flow} of $T$.
 If two group actions are orbit equivalent then their associated flows are isomorphic.
 The associated flow is transitive and free if and only if $T$ admits an  $\sigma$-finite invariant $\mu$-equivalent measure.
 In this case $T$ is said to be {\it of type $II$}.
 If the invariant measure is finite, $T$ is said to be of type $II_1$.
 If the invariant measure is infinite, $T$ is said to be of type $II_\infty$.
 If $T$ does not admit an invariant equivalent measure then $T$ is said to be {\it of type $III$}.
 Type $III$ admits further classification into subtypes $III_\lambda$, $0\le\lambda\le 1$.
 If the associated flow of $T$ is periodic with period $-\log\lambda$ for some $\lambda\in(0,1)$ then $T$ is said to be {\it of type $III_\lambda$}.
 If the associated flow is trivial (on a singletone) then $T$ is said to be {\it of type $III_1$}.
 Equivalently, $T$ is of type $III_1$ if and only if the Radon-Nikodym cocycle of $T$ is ergodic.
 If $T$ is of type $III$ but not of type $III_\lambda$ for any $\lambda\in(0,1]$ then $T$ is said to be {\it of type $III_0$.}

 \proclaim{Lemma 1.1} Let $\mu$ be invariant under $T$.
  Let  $H$ be discrete and countable. 
  Let $\alpha:\Cal R\to H$ be an ergodic cocycle.
  Then the following holds.
 \roster
 \item"\rom(i)" 
 The subrelation $\Cal R_0:=\{(x,y)\in\Cal R\mid \alpha(x,y)=1\}$ of $\Cal R$ is ergodic.
 \item"\rom(ii)" 
 For each $h\in H$, there is an element $r_h\in[\Cal R]$ such that $\alpha(r_hx,x)=h$ for $\mu$-a.e. $x\in X$.
 \endroster
 \endproclaim
 \demo{Idea of the proof}
 Pass to the ergodic skew product extension $T(\alpha)$ and use the following Hopf lemma:
 if $D=(D_h)_{h\in H}$ is an ergodic $H$-action of type $II$ and  $\lambda$ is a $D$-invariant equivalent measure then for all subsets $A$ and $B$  with $\lambda(A)=\lambda(B)$, there is a Borel bijection $\tau:A\to B$  such that the graph of $\tau$ is a subset of  the $D$-orbit equivalence ralation.
 \qed
 \enddemo
 
   Let $S$ be a nonsingular $H$-action on a standard probability space $(Y,\goth Y,\nu)$.
   Given a cocycle $\alpha:\Cal R\to H$, we can form a {\it skew product action} $T(\alpha,S)=(T(\alpha,S)_g)_{g\in G}$ of $G$ on the product space $(X\times Y,\mu\times\nu)$ by setting
   $$
   T(\alpha,S)_g(x,y):=(T_gx,S_{\alpha(T_gx,x)}y).
   $$
  Then $T(\alpha,S)$ is $(\mu\times\nu)$-nonsingular.

 \proclaim{Lemma 1.2} Let $T,\mu,H$ be as in Lemma~1.1.
 If $\alpha$ is ergodic and  $S$ is ergodic then $T(\alpha,S)$ is also ergodic. 
The associated flow of $T(\alpha,S)$ is isomorphic to the
 associated flow of $S$.
 In particular, the type of $T(\alpha,S)$ equals the type of $S$.
 \endproclaim
 \demo{Proof}
 \comment 
 Take  $\lambda\in \Cal P_G(X\times Y,\mu)$.
 Consider the disintegration $\lambda=\int_X\delta_x\times\lambda_x\,d\mu(x)$ of $\lambda$ with respect to $\mu$.
 Then $\lambda_{T_gx}=\lambda_x\circ S_{\alpha(T_gx,x)}^{-1}$ for a.e. $x\in X$.
 Since $\lambda_{x'}=\lambda_x$ for all $(x,x')\in\Cal R_0$, it follows from Lemma~1.1(i) that there is $\xi\in\Cal P(Y)$ with $\lambda_x=\xi$ for a.e. $x$.
 By Lemma~1.1(2), $\xi\in\Cal P_{H,0}(Y)$.
 \endcomment
 Let $F:X\times Y\to\Bbb R$ be a  Borel function.
 If $F$ is $T(\alpha,S)$-invariant then $F(x,y)=F(x',y)$ if $(x',x)\in\Cal R_0$ for a.e. $y$.
 By Lemma~1.1(i), $F(x,y)=f(y)$ for some Borel  function $f:Y\to \Bbb R$.
 Lemma~1.1(ii) now yields that $f$ is invariant under $S$.
 Since $S$ is ergodic,  $f$ is constant mod  $\nu$.
 Thus $F$ is constant mod $\mu\times\nu$.
 Hence $T(\alpha,S)$ is ergodic.
 
 The second claim of the lemma follows from Lemma~1.1,  the fact that $\Cal R$ is generated by $\Cal R_0$ and the family of transformations $(r_h)_{h\in H}$ and that $r_h[\Cal R_0]r_h^{-1}=[\Cal R_0]$ for each $h\in H$.
\qed
 \enddemo

 We now  recall the definition of  strongly ergodic actions (see \cite{CoWe}, \cite{JoSc} and references therein).
 Let $T$ be an ergodic  nonsingular $G$-action on  non-atomic probability space $(X,\goth B,\mu)$.
 A sequence $(B_n)_{n\in\Bbb N}$ in $\goth B$ is called {\it asymptotically invariant} if $\lim_{n\to\infty}\mu(B_n\triangle T_gB_n)=0$ for every $g\in G$.
 If every asymptotically invariant sequence $(B_n)_{n\in\Bbb N}$ is trivial, i.e. $\lim_{n\to\infty}\mu(B_n)(1-\mu(B_n))=0$ then $T$ is called {\it strongly ergodic}.  We note the the strong ergodicity is invariant under the orbit equivalence.
 We will need the following lemma.
 
  \proclaim{Lemma 1.3} 
 \roster
 \item"\rom(i)"
 If $G$ does not have property (T) then there is an ergodic  probability preserving free action $T$ of $G$ which is not strongly ergodic  \cite{CoWe}.
 \item"\rom(ii)"
  If  $T$ is an ergodic nonsingular free action of $G$ which is not strongly ergodic then for each countable discrete Abelian group $A$,  there is an ergodic  cocycle of the $T$-orbit equivalence relation with values in $A$
 (see  \cite{Sc, Corollary 1.5} and Theorem~A2 below\footnote{Since the proof of  \cite{Sc, Corollary 1.5} was not completed there we provide a complete proof of it in Appendix~A.}).
 \endroster
 \endproclaim

 \comment
 
 We now recall the definition of  AT-flows \cite{CoWo}.
Let $(n_j)_{j\in\Bbb N}$ be an infinite sequence of integers greater than $1$.
Let $X=\bigotimes_{j\in\Bbb N} \Bbb Z/n_j\Bbb Z$.
We furnish $X$ with the standard product topology.
Then $X$ is a compact Abelian group under the pointwise addition.
The countable group $G=\bigoplus_{j\in\Bbb N}\Bbb Z/n_j\Bbb Z$ embeds naturally in $X$ as a
dense subgroup.
Denote by $S$ the action of $G$ on $X$ by rotations.
Let $\nu=\bigotimes_{j\in\Bbb N}\nu_j$, where $\nu_j$ is a probability on $\Bbb Z/2\Bbb Z$.
Suppose that $\nu$ has no atoms.
We note that $S$ is $\nu$-nonsingular and ergodic.
An ergodic nonsingular $\Bbb R_+^*$-action is called  an {\it AT-flow} if it is conjugate to
the associated flow of $S$  for some choice of $(n_j)_{j\in\Bbb N}$ and $(\nu_j)_{j\in\Bbb N}$.

 \endcomment

 \head 2. Proof of the main result
 \endhead
 
The following lemma is almost a literal repetition of \cite{Bo--Ta, Lemma~4.1}, where it was proved under an additional   assumption that $T$ is measure preserving.

 \proclaim{Lemma 2.1 (Entropy addition formula)} Let $\kappa$ be a probability on $G$ and let $T$ be a nonsingular action of $G$ on a standard probability space $(X,\goth B,\mu)$.
 Given a discrete countable group $H$ and a nonsingular action $S=(S_h)_{h\in H}$ of  $H$
 on a standard probability space $(Y,\goth F,\nu)$, 
  let $\kappa_x$ denote the pushforward of $\kappa$ under the map $G\ni g\mapsto\alpha(T_gx,x)\in H$ for each  $x\in X$.
 Then 
 $$
  h_{\kappa}(T(\alpha),\mu\times\nu)=h_\kappa(T,\mu)+ \int_Xh_{\kappa_x}(S,\nu)d\mu(x).
  $$
 \endproclaim
 \demo{Proof}
 $$
 \align
 h_{\kappa}(T(\alpha,S),\mu\times\nu)&=-\sum_{g\in G}\kappa(g)\int_{X\times Y}\log\left(\frac{d(\mu\times\nu)\circ T_g(\alpha)}{d(\mu\times\nu)}(x,y)\right)d\mu(x) d\nu(y)\\
 &=h_\kappa(T,\mu)-
\int_X \sum_{g\in G}\kappa(g)\int_Y
 \log\left(\frac{d\nu\circ S_{\alpha(T_gx,x)}}{d \nu}(y)\right)d\nu(y)d\mu(x)\\
 &=h_\kappa(T,\mu)-
\int_X \sum_{h\in H}\kappa_x(h)\int_Y
 \log\left(\frac{d\nu\circ S_{h}}{d \nu}(y)\right)d\nu(y)d\mu(x)\\
 &=h_\kappa(T,\mu)+ \int_Xh_{\kappa_x}(S,\nu)d\mu(x).\qed
 \endalign
 $$
 \enddemo

 \comment
 \proclaim{Theorem 1.5}
  Let $\kappa$ be a generating measure on $G$.
  \roster
  \item
  Then for each $\epsilon>0$ and $\lambda\in [0,1]$, there is an ergodic nonsingular action $T$ of  $G$ of type $III_\lambda$ such that $0<h_\kappa(T)<\epsilon$.
  \item
  More generally,
  let $T$ be an ergodic nonsingular action of $G$ on a non-atomic probability space.
  Then for each $\epsilon>0$ and, there is a type $III$ ergodic nonsingular extension $\widehat T$ of  $T$ such that  $h_\kappa(T)<h_\kappa(\widehat T)<h_\kappa(T)+\epsilon$.
  \endroster
 \endproclaim
 \endcomment
 
 \demo{Proof of Main Theorem}
 We will proceed in two steps.
 On the first step, for each  $\epsilon>0$, we construct an ergodic nonsingular $G$-action of type $III_1$  (or of type $III_\lambda$ for some $\lambda\in(0,1)$) whose  $\kappa$-entropy is less then $\epsilon$.
 On the second step we show how to change the
 quasiinvatiant measure for the action constructed on the first step with appropriate equivalent measures
 to forse
the $\kappa$-entropy to attain all the values from the interval $(\epsilon,+\infty)$.   
%Since the isomorphism class of the associated flow is invariant under such a change, the statement of the theorem follows. 
 
{\it Step 1.} Fix $\epsilon>0$ and $\lambda\in(0,1)$.
 Fix an enumeration 
 $G=\{g_n\mid n\in\Bbb N\}$ and   a sequence of integers $1=l_1\le l_2\le\cdots$ such that $l_{n+1}-l_n\le 1$ for all $n\in\Bbb N$,   $l_n\to\infty$ and
  $$
  \sum_{n=1}^\infty\kappa(g_n)(l_n+1)<2.\tag2-1
$$
 By Lemma~1.3(i),  there is a measure preserving free action $T$ of $G$ on a standard probability space $(X,\goth B,\mu)$ which is not strongly ergodic.
 Denote by $\Cal R$ the $T$-orbit equivalence relation.
 Let $F:=\bigoplus_{n\in\Bbb N}\Bbb Z/2\Bbb Z$.
 We consider the  elements of $F$ as $\Bbb Z/2\Bbb Z$-valued functions on $\Bbb N$ with finite support.
   Given $f\in F$, we let 
   $$
   \|f\|:=\max\{j\in\Bbb N\mid f(j)\ne 0\}.
   $$
   By Lemma~1.3(ii), there exists an ergodic   cocycle $\alpha:\Cal R\to F$.
For each $n\in\Bbb N$, we can choose  $M_n\in\Bbb N$ such that
$$
\mu\bigg(\bigg\{x\in X\,\bigg|\,\max_{1\le i\le n} \|\alpha(T_{g_i}x,x)\|< M_n\bigg\}\bigg)>1-\frac 1{n2^n}.\tag2-2
$$
Without loss of generality we may assume that 
$M_n=M_{n+1}$
if and only if $l_n=l_{n+1}$ for each $n\in\Bbb N$.
Let 
$
N:=\{f\in F\mid f(M_n)=0\text{ for each  }n\in\Bbb N\}.
$
Then $N$ is a subgroup of $F$.
The quotient group  $F/N$ is identified naturally with the ``complimentary  to $F$'' subgroup $\{f\in F\mid f(n)=0\text{ for each  }n\neq M_1,M_2,\dots\}$ which is, in turn,  isomorphic  to $F$ in a natural way.
Hence passing from $\alpha$ to the quotient cocycle 
$$
\alpha+ N:\Cal R\ni(x,y)\mapsto\alpha(x,y)+N\in F/N
$$
 means that we may assume without loss of generality that $M_n=l_n$  for each $n\in\Bbb N$ 
 in~\thetag{2-2}.
(We use here a simple fact that $\alpha+ N$ is ergodic whenever  $\alpha$ is.)
Therefore applying \thetag{2-2} we  obtain that 
$$
\aligned
\int_X\|\alpha(T_{g_{n}}x,x)\|\,d\mu(x)&=\sum_{s=1}^\infty s \mu(\{x\in X\mid \|\alpha(T_{g_{n}}x,x)\|=s\})\\
&\le l_n+\sum_{s>l_n}s \mu(\{x\in X\mid \|\alpha(T_{g_{n}}x,x)\|=s\})\\
&\le
l_n+\sum_{s>l_n}\frac 1{2^{s}}\\
&\le l_n+1.
\endaligned
\tag2-3
$$
 \comment
 We find $L>0$ such that the subset
 $$
X_{L}:=\Bigg\{x\in X\,\Bigg| \sum_{\{h\mid N_h\not\subset [1,L]\}}\xi_x(h)<\epsilon\Bigg\}
 $$
 has $\mu$-measure greater then $1-\epsilon$
 \endcomment
 The second inequality here  follows from the fact that for each $s>l_n$, we have $s=l_m$ for some $m\ge n$ and hence
 $$
  \mu(\{x\in X\mid \|\alpha(T_{g_{n}}x,x)\|=s\})\le \frac1{m2^m}\le \frac1{s2^{s}}.
  $$
  because $m\ge l_m=s$.
  
Now we consider  $F$ as a (dense) subgroup  of the compact Abelian group $K:=(\Bbb Z/2\Bbb Z)^\Bbb N$ of all $\Bbb Z/2\Bbb Z$-valued functions on $\Bbb N$.
 Denote by $S$ the action of $H$ on $K$ by translations.
 Let $\nu_n$ denote the  distribution on $Z/2\Bbb Z$ such that 
 $$
 \nu_n(0)=\frac 1{1+ e^{\epsilon_n} },\quad   \nu_n(1)=\frac { e^{\epsilon_n}}{1+e^{\epsilon_n} } 
 $$
  for some sequence $(\epsilon_n)_{n\in\Bbb N}$ of reals 
 such that  $\lim_{n\to\infty}\epsilon_n=0$, $\sum_{n\in\Bbb N}\epsilon_n^2=\infty$ and
 $\max\{|\epsilon_n|\mid n\in\Bbb N\}<\epsilon$.
\comment
 $$
\delta:= \sup_{n\in\Bbb N}\left\{(\nu_n(1)-\nu_n(0)) \log\frac{\nu_n(1)}{\nu_n(0)}\right\}<\frac \epsilon 2.
 $$
 \endcomment
 Let  $\nu=\bigotimes_{n\in\Bbb N}\nu_n$.
 Then $S$ is  $\nu$-nonsingular,  ergodic
 and of type $III_1$
 \cite{ArWo}.
 Therefore by Lemma~1.2, the skew product $G$-action $T(\alpha,S)$ on $(X\times K,\mu\times\nu)$ is ergodic and of type $III_1$.
 Hence $h_\kappa(T(\alpha,S))>0$.
 To estimate $h_\kappa(T(\alpha,S))$ from above, we first let $N_f:=\{n\in\Bbb N\mid f(n)\ne 0\}$ for  $f\in F$.
 It is obvious that $\# N_f\le\|f\|$.
 %Without loss of generality we may assume that $1$ is not a limit point of the sequence $(\nu_n(1)/%\nu_n(0))_{n\in\Bbb N}$.
Since
 $$
 \align
 -\int_Y\log\left(\frac{d\nu\circ S_f}{d\nu}\right)(y)d\nu(y)&=-\sum_{n\in N_f}\int_{\Bbb Z/2\Bbb Z}\log\left(\frac{\nu_n(y_n+1)}{\nu_n(y_n)}\right)d\nu_n(y_n)\\
 &=
\sum_{n\in N_f}(\nu_n(1)-\nu_n(0)) \log\frac{\nu_n(1)}{\nu_n(0)},
\endalign
 $$
 we obtain that  for each probability $\xi$ on $F$,
 $$
 \aligned
 h_\xi(S,\nu)&=\sum_{f\in F}\xi(f)\sum_{n\in N_f}(\nu_n(1)-\nu_n(0)) \log\frac{\nu_n(1)}{\nu_n(0)}\\
&\le\sum_{f\in F}\xi(f)\sum_{n\in N_f}|\epsilon_n|\\
 & \le\epsilon\|f\|\sum_{f\in F}\xi(f).
 \endaligned
 \tag2-4
 $$
 Since $T$ preserves $\kappa$, it follows that $h_\kappa(T,\mu)=0$.
 Then Lemma~2.1, \thetag{2-4}  and \thetag{2-3}  yield that
 $$
 \align
 h_{\kappa}(T(\alpha,S),\mu\times\nu)&\le\epsilon\int_X\sum_{f\in F}\kappa_x(f)\|f\|\,d\mu(x)\\
 %&\le\sup_{n\in\Bbb N}d_{\text{KL}}(\nu_n^*|\nu_n)
 %\int_X\sum_{f\in F}\kappa_x(f)\# N_f\,d\mu(x)\\
&=\epsilon\sum_{g\in G}\kappa(g)
 \int_X\|\alpha(T_gx,x)\|\,d\mu(x)\\
 &\le  \epsilon\sum_{n=1}^\infty\kappa(g_n)(l_n+1).
 \endalign
 $$
 It now follows from \thetag{2-1} that
  $
  h_{\kappa}(T(\alpha,S),\mu\times\nu)\le2\epsilon.
 $
Hence
 $$
 \inf\{h_\kappa(A)\mid \text{$A$ is a type $III_1$ ergodic free action of $G$}\}=0.
 $$
 In a similar way we may show that
 $$
 \inf\{h_\kappa(A)\mid \text{$A$ is a type $III_\lambda$ ergodic free action of $G$, $\lambda\in(0,1)$}\}=0.
 \tag2-5
 $$ 
 For that we argue as above but with a different measure $\nu$.
 Indeed,
 let $\nu_n$ denote the  distribution on $Z/2\Bbb Z$ such that 
 $$
 \nu_n(0)=\frac 1{1+ e^{\epsilon}} ,\quad   \nu_n(1)=\frac {e^{\epsilon}}{1+e^{\epsilon}} 
 $$
% that 
  %$
%\delta:= \sup_{n\in\Bbb N}\left\{(\nu_n(1)-\nu_n(0)) \log\frac{\nu_n(1)}{\nu_n(0)}\right\}<\frac \epsilon 2.
 %$
  Let  $\nu=\bigotimes_{n\in\Bbb N}\nu_n$.
 Then $S$ is  $\nu$-nonsingular,  ergodic
 and of type $III_{e^{-\epsilon}}$
 \cite{ArWo}.
 Therefore by Lemma~1.2, the skew product $G$-action $T(\alpha,S)$ on $(X\times K,\mu\times\nu)$ is ergodic and of type $III_{e^{-\epsilon}}$.
 Hence $h_\kappa(T(\alpha,S))>0$.
As in the $III_1$-case considered above, we obtain that
$
  h_{\kappa}(T(\alpha,S),\mu\times\nu)<2\epsilon
 $
 and hence \thetag{2-5} follows.

 {\it Step 2.} Given $\epsilon>0$, let $\nu$ be a measure on $K$ such that 
 $$
 h_\kappa(T(\alpha,S),\mu\times\nu)<\epsilon.\tag2-5
 $$
  We choose $n_0>0$ such that
  $$
  \int_X\kappa_x(\{f\in F\mid f(n_0)\ne 0\})\,d\mu(x)>0.\tag2-6
  $$
  It exists because otherwise we would have that $\kappa_x$ is supported at $0$  for a.e. $x\in X$.
  The latter yields that $\alpha(T_gx,x)=0$ at a.e. $x$ for all $g$ from the support of $\kappa$.
  Since $\kappa$ is generating, it follows that $\alpha$ is trivial, a contradiction.
  
  Let $\omega$ be a probability on $\Bbb Z/2\Bbb Z$ supported at   $1$\footnote{We consider the group $\Bbb Z/2\Bbb Z$ as $\{0,1\}$ with addition mod 2.}.
  For each $\theta\in(0,1]$, we let
  $$
  \nu_n^\theta=\cases
  \nu_n, &\text{if }n\ne n_0\\
  \theta\nu_n+(1-\theta)\omega,
  &\text{if }n= n_0
  \endcases
  $$
  and $\nu^\theta:=\bigotimes_{n\in\Bbb N}\nu^\theta_n$. 
  Then  $\nu^\theta$ is equivalent to $\nu$ and hence  $\mu\times\nu^\theta$ is equivalent to $\mu\times\nu$. 
  Therefore  the dynamical systems  $(T(\alpha,S),\mu\times\nu^\theta)$  and  $(T(\alpha,S),\mu\times\nu)$ are of the same Krieger's type.
  It follows from the equality in \thetag{2-4} that
  $$
  h_{\kappa_x}(S,\nu)-  h_{\kappa_x}(S,\nu^\theta)=
\kappa_x(\{f\in F\mid f(n_0)\ne 0\}) (\Phi(\nu_{n_0}(0))-\Phi(\nu_{n_0}^\theta(0))),
   $$
  where
  $\Phi(t):=(1-2t)\log\frac{1-t}{t}$, if $t\in(0,1)$.
Therefore the map 
  $$
  (0,1]\ni \theta\mapsto h_{\kappa}(T(\alpha,S),\mu\times\nu^\theta)= \int_Xh_{\kappa_x}(S,\nu^\theta)\,d\mu(x)\in\Bbb R
  $$
  is continuous.
  In view of \thetag{2-6},
 this map goes to infinity as $\theta\to 0$.
Since $\nu^1=\nu$ and \thetag{2-5} holds,
  it follows that  
  $$
  \{h_\kappa(T(\alpha,S), \mu\times\nu^\theta)\mid\theta\in(0,1]\}\supset(\epsilon,\infty),
  $$
  as desired.
  \qed
  \enddemo

  \head Appendix A
  \endhead

  \comment
  \proclaim{Proposition}
  Let $\Cal R$ be an ergodic measure preserving countable equivalence relation on a standard probability space $(X,\goth B,\mu)$.
  Let $F=\bigoplus_{i\in\Bbb N}G_i$, where each $G_i$ is $\Bbb Z/2\Bbb Z$.
  If there is a cocycle $\alpha:\Cal R\to F$ which is not cohomologous to a cocycle of $\Cal R$ with values in a finite subgroup of $F$ then there is an ergodic cocycle $\beta:\Cal R\to F$.
  \endproclaim
  
  \demo{Proof}
  For each $i\in\Bbb N$, we denote by $\alpha_i:\Cal R\to G_i$ the projection of $\alpha$ to the $i$-th coordinate.
  Of course, each cocycle $\alpha_i$ is either ergodic or coboundary.
  It follows from the assumption of the lemma that the set $I_1:=\{i\in\Bbb N\mid\text{$\alpha_i$ is ergodic}\}$
  is infinite.
  Fix $i_1\in I_1$.
  For each $i\in I_1\setminus\{i_1\}$, we denote by $\alpha_{i_1,i}$ the projection of $\alpha$ to $G_{i_1}\oplus G_i$.
  Then there is a subgroup $H_{i_1,i}\subset G_{i_1}\oplus G_i$ such that  $\alpha_{i_1,i}$ is cohomologous to an ergodic cocycle with values in $H_{i_1,i}$.
  Since $i_1,i\in I_1$, it follows that the projection of $H_{i_1,i}$ to each of the two coordinates is onto, i.e. $H_{i_1,i}$ is an ``algebraic joining''  of $G_{i_1}$ and $G_i$.
  Hence either $H_{i_1,i}=\{(g,g)\mid g\in G_i\}$ or $H_{i_1,i}=G_{i_1}\times G_i$.
  The former case means that $\alpha_i$ is cohomologous to $\alpha_{i_1}$.
  We claim that the set
  $\{i\in I_1\mid i\ne i_1, H_{i_1,i}=\{(g,g)\mid g\in G_i\}\}$ is finite.
  \enddemo
  \endcomment
  
 Let $\Cal R$ be an ergodic measure preserving countable equivalence relation on a nonatomic standard probability space $(X,\goth B,\mu)$ and let $G$ be a locally compact second countable group. 
 A cocycle $\rho:\Cal R\to G$ is called {\it regular} if the action of $G$ associated with $\rho$ is transitive.
 For instance, an ergodic cocycle is regular. 
 A coboundary is also regular.

  \proclaim{Proposition A1} %Let $\Cal R$ be an ergodic measure preserving countable equivalence relation on a standard nonatomic probability space $(X,\goth B,\mu)$.
  Let $A$ be an amenable discrete countable group and let $H$ be a locally compact second countable amenable group.
  Let $\alpha:\Cal R\to A$ be a  cocycle.
  If $\alpha$ is not regular then there is an ergodic cocycle of $\Cal R$ with values in $H$.
  \endproclaim
  
  \demo{Proof}
  Let $\Cal A$ stand for the ``transitive'' equivalence relation on $A$, i.e. $\Cal A=A\times A$.
  Let $\lambda$ be a probability measure on $A$ which is equivalent to Haar measure.
  Then the equivalence relation $\Cal R\times\Cal A$ is an ergodic equivalence relation on the probability space $(X\times A,\mu\times\lambda)$. 
  Let $V=(V_a)_{a\in A}$ denote the nonsingular action of $A$ on  $(X\times A,\mu\times\lambda)$  by right rotations  along the second coordinate.
  Then $\Cal R\times\Cal A$ is generated by a subrelation $\Cal R(\alpha)$ and $V$, i.e. two points
  $z_1,z_2\in X\times A$ are $(\Cal R\times\Cal A)$-equivalent if and only if the points $V_{a_1}z_1$ and $V_{a_2}z_2$
are $\Cal R(\alpha)$-equivalent for some $a_1,a_2\in A$.
Let $W=(W_a)_{a\in A}$ stand for the action of $A$ associated with $\alpha$.
Denote by $(\Omega,\nu)$ the space of this action.
Then we can assume that there is a Borel  map $\pi:X\times A\to \Omega$ 
such that $\nu=(\mu\times\lambda)\circ \pi^{-1}$ and $\pi\circ V_a=W_a\circ \pi$ for each $a\in A$.
We observe that $\pi$ is the $\Cal R(\alpha)$-ergodic decomposition of $X\times A$. 
Since $A$ is amenable, the $W$-orbit equivalence relation  $\Cal I$ on $(\Omega,\nu)$ is hyperfinite \cite{Co--We}.
It is ergodic.
Since $\alpha$ is not regular, $\Cal I$ is non-transitive.
Hence there is an ergodic cocycle $\beta:\Cal I\to H$ (see \cite{He} and \cite{GoSi}).
We now define a cocycle $\beta^*:\Cal R\times \Cal A\to H$ by setting
$$
\beta^*(z_1,z_2):=\beta(\pi( z_1),\pi( z_2)).
$$
Then $\beta^*$ is well defined.
Since  $\pi$ is the $\Cal R(\alpha)$-ergodic decomposition, it follows that $\beta^*$ is ergodic.
Restricting $\Cal R\times \Cal A$ and $\beta^*$ to the subset $X\times\{1_A\}$ we obtain $\Cal R$ and a cocycle of $\Cal R$
with values in $H$ respectively.
Of course, this cocycle is also ergodic.
\qed
  \enddemo

  Let $A$ be an Abelian locally compact noncompact second countable group.
  For a cocycle $\alpha:\Cal R\to A$, we denote by $E(\alpha)\subset A\sqcup\{\infty\}$ the essential range of $\alpha$ (see \cite{Sc1, Definition 3.1}).
  The following theorem provides a complete proof of  \cite{Sc2, Corollary~1.5}  (it was assumed additionally in \cite{Sc2} that  $A$ is Abelian).

  \proclaim{Theorem A2} Let $T=(T_g)_{g\in G}$ be an ergodic measure preserving action of $G$ on a standard nonatomic probability space $(X,\goth B,\mu)$.
  Let $A$ be a countable amenable group.
  If $T$ is not strongly ergodic then there is an ergodic  cocycle of the $T$-orbit equivalence relation $\Cal R$ with values in $A$. 
  \endproclaim

  \demo{Proof}
  Let $\Cal R$ denote the $T$-orbit equivalence relation.
  Let $K$, $F$, $S$, $N_f$ and $\|.\|$ denote the same objects as in the proof of Main Theorem.
  Let $\lambda_K$ stand for the Haar measure on $K$.
 For each $n\in\Bbb N$, denote by $f_n$ the element of $F$ such that $N_{f_n}=\{n\}$.
  By \cite{JoSc, Lemma~2.4}, there are
  a Borel map  $\pi:X\to K$ and
   a sequence $(V_n)_{n\in\Bbb N}$ of transformations in $[\Cal R]$
   such that $\mu\circ \pi^{-1}=\lambda_K$, 
   $$
   \{\pi(T_gx)\mid g\in G\}=\{S_f\pi(x)\mid f\in F\}
   $$
   and
    $\pi(V_nx)=S_{f_n}\pi(x)$ for a.a. $x\in X$.
 Denote by $\Cal S$ the $S$-orbit equivalence relation on $K$.
  We define a cocycle $\beta:\Cal S\to F$ by setting
  $\beta(S_fy,y):=f$, $y\in Y$, $f\in F$.
 By \cite{Sc1, Proposition~3.15}, $E(\beta)=\{0,+\infty\}$.
We now define a cocycle $\beta^*:\Cal R\to F$ by setting
  $$
  \beta^*(T_gx,x):=\beta(\pi(T_gx),\pi(x)),
  $$
   $x\in X$, $g\in G$.
   Since $E(\beta^*)\subset E(\beta)$, we obtain  that either $E(\beta^*)=\{0,\infty\}$ or $E(\beta^*)=\{0\}$.
  In the latter case, $\beta^*$ is a coboundary \cite{Sc1}.
  Hence there is a Borel map $\xi:X\to K$ such that $\beta^*(V_nx,x)=\xi(V_nx)-\xi(x)$ for each $n$ at a.e. $x\in X$.
  There is a constant $C>0$ and a subset $X_C\subset X$ such that $\mu(X_C)>3/4$ and
  $\|\xi(x)\|<C$ for all $x\in X_C$.
  It follows that for each $n$, 
  $$
 \sup_{x\in X_C\cap V_n^{-1}X_C} \|\beta^*(V_nx,x)\|<C.
  $$
This contradicts to the fact that $\|\beta^*(V_nx,x)\|=n$ for a.a. $x\in X$ and $n\in\Bbb N$.
Hence $E(\beta^*)=\{0,\infty\}$.
This yields that $\beta^*$ is not regular.
It remains to apply~Proposition~A1.
\qed
  \enddemo

  \comment
  Since $V$ normalizes $\Cal R(\alpha)$, i.e. two points $z_1,z_2\in X\times A$ are $\Cal R(\alpha)$ equivalent if and only if the points $V(a)z_1$ and $V_a(z_2)$ are $\Cal R(\alpha)$-equivalent for each $a\in A$,
  \endcomment

\enddocument